\documentclass{amsart}
\usepackage{amssymb}
\usepackage{mathrsfs}

\newcommand{\bigpart}[1]{\par\bigskip\bigskip\goodbreak\noindent{\large\textbf{#1}}\par\medskip\noindent}
\newcommand{\gp}{\mbox{-\textit{\tiny gp}}}
\newcommand{\gpbl}{{\mbox{\textit{\tiny gp}}}}

\long\def\forget#1\forgotten{}

\newcommand{\scrA}{\mathscr{A}}
\newcommand{\scrB}{\mathscr{B}}

\newcommand{\MF}{{M\!\!F}}
\newcommand{\cC}{\mathcal{C}}
\newcommand{\cU}{\mathcal{U}}
\newcommand{\cV}{\mathcal{V}}

\newcommand{\cP}{\mathcal{P}}

\newcommand{\seq}[1]{\{#1\}_{n\in\N}}

\newcommand{\Cal}{\mathcal}

\newcommand{\Tau}{\mathrm{T}}

\newcommand{\cF}{{\Cal F}}

\newcommand{\M}{\mathcal{M}}
\newcommand{\N}{\mathbb{N}}

\newcommand{\NN}{{{}^\N\N}}

\renewcommand{\O}{\Cal O}

\newcommand{\Union}{\bigcup}

\newcommand{\Impl}{\Rightarrow}

\renewcommand{\b}{{\mathfrak b}}

\renewcommand{\c}{{\mathfrak c}}
\renewcommand{\d}{{\mathfrak d}}

\renewcommand{\i}{\item}

\newcommand{\oo}{\infty}

\newcommand{\Iff}{\Leftrightarrow}

\newcommand{\nin}{\not\in}
\newcommand{\cat}{\hat{\ }}
\newcommand{\sbst}{\subseteq}
\newcommand{\spst}{\supseteq}
\newcommand{\sm}{\setminus}

\newcommand{\as}{\sbst^*}
\renewcommand{\|}{\upharpoonright}

\newcommand{\<}{\langle}
\renewcommand{\>}{\rangle}

\newcommand{\cov}{{\sf cov}}
\newcommand{\add}{{\sf add}}

\newtheorem{thm}{Theorem}[section]
\newtheorem{prop}[thm]{Proposition}
\newtheorem{prob}[thm]{Problem}

\newtheorem{lem}[thm]{Lemma}
\newtheorem{cor}[thm]{Corollary}
\newtheorem{conj}{Conjecture}

\theoremstyle{definition}
\newtheorem{definition}[thm]{Definition}

\theoremstyle{remark}
\newtheorem{rem}[thm]{Remark}

\newcommand{\be}{\begin{enumerate}}
\newcommand{\ee}{\end{enumerate}}
\newcommand{\bi}{\begin{itemize}}
\newcommand{\ei}{\end{itemize}}

\newcommand{\sone}{{\sf S}_1}    \newcommand{\sfin}{{\sf S}_{fin}}
\newcommand{\ufin}{{\sf U}_{fin}}

\newcommand{\gone}{{\sf G}_1}    \newcommand{\gfin}{{\sf G}_{fin}}

\newcommand{\lft}[2]{\mathopen\ifcase#1{}\oo\or
                        \big#2\or\Big#2\else\oo\fi}
\newcommand{\rgt}[2]{\mathclose\ifcase#1{}\oo\or
                        \big#2\or\Big#2\else\oo\fi}

\title{Strong $\gamma$-sets and other singular spaces}
\author{Boaz Tsaban}
\thanks{This paper constitutes a part of the author's doctoral dissertation at the Bar-Ilan University.}
\address{Department of Mathematics,
Weizmann Institute of Science, Rehovot 76100, Israel}
\email{boaz.tsaban@weizmann.ac.il}
\urladdr{http://www.cs.biu.ac.il/\~{}tsaban}



\begin{document}
\begin{abstract}
Whereas
the Gerlits-Nagy $\gamma$ property
is strictly weaker than
the Galvin-Miller strong $\gamma$ property,
the corresponding
strong notions for the
Menger, Hurewicz, Rothberger, Gerlits-Nagy $(*)$,
Arkhangel'ski\v{i} and Sakai
properties are equivalent to the original ones.
The main result is that almost each of these properties
admits the game theoretic characterization
suggested by the stronger notion.
We also solve a related problem of Ko\v{c}inac and Scheepers,
and answer a question of Iliadis.
\end{abstract}

\maketitle

\section{Introduction}

\subsection{Thick covers}
Let $X$ be an infinite topological space.
Throughout this paper, by \emph{open cover}
we mean a \emph{countable}
collection $\cU$ of open subsets of $X$ such that
$\cup\cU=X$ and
$X\nin\cU$.
The focus on countable covers
allows us to have no restrictions at all on the topology of $X$.\footnote{A
standard alternative approach is to consider spaces $X$ such that all
finite powers of $X$ are Lindel\"of. This guarantees that each cover
of a type considered here contains a countable cover of the same type.
}
This is useful when we wish to project the results into the purely combinatorial case --
Section \ref{purecomb}.
The additional restriction that $X\nin\cU$ is to avoid trivialities.

Let $\cU$ be an open cover of $X$.
$\cU$ is an \emph{$n$-cover} of $X$ if for each $F\sbst X$ with $|F|\le n$,
there is $U\in\cU$ such that $F\sbst U$.
$\cU$ is an \emph{$\omega$-cover} of $X$ if for each finite $F\sbst X$,
there is $U\in\cU$ such that $F\sbst U$ (that is, $\cU$ is an $n$-cover of $X$
for each $n$).
$\cU$ is a \emph{$\gamma$-cover} of $X$ if
each element of $X$ belongs to all but finitely many
members of $\cU$.

\subsection{$\gamma$-sets and strong $\gamma$-sets}
According to Gerlits and Nagy \cite{GN}, a topological space $X$ is a
\emph{$\gamma$-set}
if each $\omega$-cover of $X$ contains a $\gamma$-cover of $X$.
Gerlits and Nagy consider the following seemingly stronger property:
\begin{quote}
For each sequence $\seq{\cU_n}$ of $\omega$-covers of $X$
there exist members $U_n\in\cU_n$, $n\in\N$, such that $\seq{U_n}$
is a $\gamma$-cover of $X$.
\end{quote}
Using Scheepers' notation \cite{coc1},
this property is a particular instance of the following selection hypothesis
(where $\scrA$ and $\scrB$ are any collections of covers of $X$):
\bi
\item[$\sone(\scrA,\scrB)$:]
For each sequence $\seq{\cU_n}$ of members of $\scrA$,
there exist members $U_n\in\cU_n$, $n\in\N$, such that $\seq{U_n}\in\scrB$.
\ei
Let $\Omega$ and $\Gamma$ denote the
collections of open $\omega$-covers and $\gamma$-covers of $X$,
respectively. Then the property considered by Gerlits and Nagy is
$\sone(\Omega,\Gamma)$, who proved that $X$ is a $\gamma$-set if,
and only if, $X$ satisfies $\sone(\Omega,\Gamma)$ \cite{GN}.

This result motivates the following definition.
According to Galvin and Miller \cite{GM}, a space $X$ is a
\emph{strong $\gamma$-set} if there exists an increasing sequence
$\seq{k_n}$ such that:
\begin{quote}
For each sequence $\seq{\cU_n}$ where for each $n$ $\cU_n$ is an open $k_n$-cover of $X$,
there exist members $U_n\in\cU_n$, $n\in\N$, such that $\seq{U_n}$ is a $\gamma$-cover of $X$.
\end{quote}
Clearly every strong $\gamma$-set is a $\gamma$-set; however the
properties are not provably equivalent (e.g., in \cite{Brendle} it is
shown that assuming CH, there exists an uncountable $\gamma$-set $X$ such
that no uncountable subset of $X$ is a strong $\gamma$-set).

As in the case of $\gamma$-sets, it will be convenient to introduce
the following general notation.
\begin{definition}
Assume that $\scrA_n$, $n\in\N$, and $\scrB$, are collections of covers of a space $X$.
Define the following selection hypothesis.
\begin{itemize}
\item[$\sone(\seq{\scrA_n},\scrB)$:]
For each sequence $\seq{\cU_n}$ where $\cU_n\in\scrA_n$ for all $n$,
there exist members $U_n\in\cU_n$, $n\in\N$, such that $\seq{U_n}\in\scrB$.
\end{itemize}
\end{definition}
For each $n$ denote by $\O_n$ the collection of all
open $n$-covers of a space $X$. Then $X$ is a strong $\gamma$-set if, and
only if, there exists an increasing sequence $\seq{k_n}$ such that
$X$ satisfies $\sone(\seq{\O_{k_n}},\Gamma)$.
This does not fit exactly into the family of properties of the form
$\sone(\seq{\scrA_n},\scrB)$, because of the external quantifier.
However, in Section \ref{Qelim} we show that this quantifier can be eliminated,
so that $X$ is a strong $\gamma$-set if, and
only if, $X$ satisfies $\sone(\seq{\O_n},\Gamma)$.
This motivates the study of the generalized selection hypothesis,
which is the aim of this paper.

The first part of the paper deals with the classical selection operators.
In Section 2, as said above, we prove quantifier elimination for the $\gamma$-property.
In Section 3 we introduce two mild assumptions on thick covers which
allow this sort of quantifier elimination.
In Sections 4 and 5 we supply a variety of examples, showing that many properties
which appear in the literature are equivalent to their stronger version.
In Section 6 we answer a question of Iliadis by showing that no new property
is obtained by considering the generalized selection hypothesis for the standard types
of covers, except for the strong $\gamma$-property.

The second part of the paper deals with game theoretic versions of the studied
properties.
In Sections 7--9 we supply new methods of reductions between game strategies,
and give new game theoretic characterizations to most of the properties mentioned
in the first part of the paper.
In Section \ref{purecomb} we describe an application of the obtained results to the purely combinatorial
case.

\bigpart{Part I: Strong versions of the classical selection operators}
\section{Strong $\gamma$-sets and quantifier elimination}\label{Qelim}

The following theorem shows that the external quantifier in the definition
of a strong $\gamma$-set can be eliminated.

\begin{thm}\label{sgnicedef}
For each space $X$, the following are equivalent:
\be
\i $X$ is a strong $\gamma$-set, that is:\\
There exists an increasing sequence $\seq{k_n}$ such that $X$ satisfies $\sone(\seq{\O_{k_n}},\Gamma)$.
\i For each increasing sequence $\seq{k_n}$, $X$ satisfies $\sone(\seq{\O_{k_n}},\Gamma)$.
\i There exists a sequence $k_n\to\infty$, such that $X$ satisfies $\sone(\seq{\O_{k_n}},\Gamma)$.
\i For each sequence $k_n\to\infty$, $X$ satisfies $\sone(\seq{\O_{k_n}},\Gamma)$.
\i $X$ satisfies $\sone(\seq{\O_n},\Gamma)$.
\ee
\end{thm}
\begin{proof}
It is clear that $4\Impl2\Impl 5\Impl 1\Impl 3$.
It remains to show that $3\Impl 4$.

Assume that $k_n\to\infty$ such that $X$ satisfies $\sone(\seq{\O_{k_n}},\Gamma)$,
and let $\seq{\cU_n}$ be such that $\cU_n\in\O_n$ for each $n$.
\begin{lem}
For $\cU\in\O_{k_1}$ and $\cV\in\O_{k_2}$, define
$$\cU\land\cV = \{U\cap V : U\in\cU, V\in\cV\}.$$
Then $\cU\land \cV$ is an open $\min\{k_1,k_2\}$-cover of $X$ refining $\cU$ and
$\cV$. (Moreover, the operation $\land$ is associative.)\hfill\qed
\end{lem}
For each $m$ let
\begin{equation}\label{minn}
n=\min\{i : m\le k_j\mbox{ for all }j\ge i\},
\end{equation}
and set
$$\cV_m=\cU_{k_n}\land\dots\land\cU_{k_{n+1}-1}.$$
Then each $\cV_m\in\O_m$. Use $\sone(\seq{\O_{k_n}},\Gamma)$ to extract
from the sequence $\seq{\cV_{k_n}}$ elements $V_{k_n}\in\cV_{k_n}$
such that $\seq{V_{k_n}}$ is a $\gamma$-cover of $X$.
For each $m$, let $n$ be as in Equation \eqref{minn}.
As $\cV_m$ refines $\cU_k$ for all $k$ where $k_n\le k<k_{n+1}$, we can choose for each such $k$
an element $U_k\in\cU_k$ such that $V_{k_n}\sbst U_k$.
For $0\le k< k_0$ choose any $U_k\in\cU_k$ (this is a finite set so we need not
worry about it). Then $\seq{U_n}$ is a $\gamma$-cover of $X$ and for each
$n$, $U_n\in\cU_n$.
\end{proof}

We now consider the following general selection hypotheses,
the first due to Scheepers and the second being a ``strong'' version
of the first.
\begin{definition}~
\bi
\item[$\sfin(\scrA,\scrB)$\index{$\sfin(\scrA,\scrB)$}:]
For each sequence $\seq{\cU_n}$ of members of $\scrA$,
there exist finite (possibly empty)
subsets $\cF_n\sbst\cU_n$, $n\in\N$, such that $\Union_{n\in\N}\cF_n\in\scrB$.
\i[$\sfin(\seq{\scrA_n},\scrB)$:]
For each sequence $\seq{\cU_n}$ where $\cU_n\in\scrA_n$ for all $n$,
there exist finite (possibly empty)
subsets $\cF_n\sbst\cU_n$, $n\in\N$, such that $\Union_{n\in\N}\cF_n\in\scrB$.
\ei
\end{definition}

In \cite{coc2} it is proved that $\sone(\Omega,\Gamma)=\sfin(\Omega,\Gamma)$.
A natural question is whether $\sone(\seq{\O_n},\Gamma)=\sfin(\seq{\O_n},\Gamma)$.
The following characterization of the $\gamma$ property answers this question.
\begin{thm}\label{gamma}
$\sone(\Omega,\Gamma)=\sfin(\seq{\O_n},\Gamma)$.
\end{thm}
As $\gamma$-sets need not be strong $\gamma$-sets,
the properties $\sone(\seq{\O_n},\Gamma)$ and $\sfin(\seq{\O_n},\Gamma)$
are not provably equivalent.

The characterization in Theorem \ref{gamma}
can be proved in a more general setting.
\begin{thm}\label{sfinequiv}
Let $\scrB$ be any collection of open covers of $X$.
Then $\sfin(\Omega,\scrB)=\sfin(\seq{\O_n},\scrB)$.
\end{thm}
\begin{proof}
It is enough to show that $\sfin(\Omega,\scrB)$ implies $\sfin(\seq{\O_n},\scrB)$.
Assume that $\seq{\cU_n}$ is a sequence of open $n$-covers of $X$.
Let $\seq{A_n}$ be a partition of $\N$ into infinitely many infinite sets.
For each $n$ define $\cV_n=\Union_{m\in A_n}\cU_m$. Then each $\cV_n$ is
an $\omega$-cover of $X$. Apply $\sfin(\Omega,\scrB)$ to extract
finite subsets $\cF_n\sbst\cV_n$, $n\in\N$, such that $\cV=\Union_{n\in\N}\cF_n\in\scrB$.
For each $n$ and each $m\in A_n$, define $\tilde\cF_m=\cF_n\cap\cU_m$.
Then for each $n$ $\tilde\cF_n$ is a finite subset of $\cU_n$, and
$\Union_{n\in\N}\tilde\cF_n=\cV\in\scrB$.
\end{proof}

\section{Other strong properties}

As mentioned in the previous section, the strong $\gamma$ property
$\sone(\seq{\O_n},\Gamma)$ is not provably equivalent to the
usual $\gamma$ property $\sone(\Omega,\Gamma)$.
Many other properties which were studied in the literature are
equivalent to properties of the form $\sone(\Omega,\scrB)$
or $\sfin(\Omega,\scrB)$ for suitably chosen $\scrB$
\cite{LecceSurvey}. We will show that for all of these properties,
the stronger versions are equivalent to them.

We first show that as in Theorem \ref{sgnicedef},
we do not get anything new if
we consider properties of the form $\sone(\seq{\O_{k_n}},\scrB)$
and $\sfin(\seq{\O_{k_n}},\scrB)$ for general increasing sequences
$k_n$.
In the case of $\sfin$ this is an immediate corollary
of Theorem \ref{sfinequiv}.
In the case of $\sone$ we need some assumptions on
$\scrB$.
\begin{definition}\label{thickcovs}
A collection $\scrA$ of open covers of a space $X$ is
\emph{finitely thick} if:
\be
\i\label{larger}
If $\cU\in\scrA$ and for each $U\in\cU$ $\cF_U$ is a finite
nonempty
family of open sets such that for each $V\in\cF_U$,
$U\sbst V\neq X$, then $\Union_{U\in\cU}\cF_U\in\scrA$.
\i\label{finadd}
If $\cU\in\scrA$ and $\cV=\cU\cup\cF$ where $\cF$ is finite
and $X\nin\cF$, then $\cV\in\scrA$.
\ee
A collection $\scrA$ of open covers of a space $X$ is
\emph{countably thick} if for each $\cU\in\scrA$ and
each countable family $\cV$ of open subsets of $X$
such that $X\nin\cV$, $\cU\cup\cV\in\scrA$.
\end{definition}
None of these two thickness properties implies the second.
The collections $\O$, $\O_n$ (for each $n$), and $\Omega$
are both finitely and countably thick.
$\Gamma$ is finitely thick but not necessarily countably thick,
and $\Lambda$, the collection of all large covers of $X$,
is countably thick but not necessarily finitely thick.

We have the following generalization of Theorem \ref{sgnicedef}.

\begin{thm}\label{generalnicedef}
Assume that
$\scrB$ is a finitely or countably thick collection of open covers of $X$.
For each space $X$, the following are equivalent:
\be
\i There exists an increasing sequence $\seq{k_n}$ such that $X$ satisfies $\sone(\seq{\O_{k_n}},\scrB)$.
\i For each increasing sequence $\seq{k_n}$, $X$ satisfies $\sone(\seq{\O_{k_n}},\scrB)$.
\i There exists a sequence $k_n\to\infty$, such that $X$ satisfies $\sone(\seq{\O_{k_n}},\scrB)$.
\i For each sequence $k_n\to\infty$, $X$ satisfies $\sone(\seq{\O_{k_n}},\scrB)$.
\i $X$ satisfies $\sone(\seq{\O_n},\scrB)$.
\ee
\end{thm}
\begin{proof}
The case where $\scrB$ is finitely thick
is proved exactly as in Theorem \ref{sgnicedef}.
The case where $\scrB$ is countably thick follows from
Theorem \ref{soneequiv}.
\end{proof}

The fact the $\Gamma$
is not countably thick is related in a straightforward manner
to the fact that $\gamma$-sets need not be strong $\gamma$-sets.
\begin{thm}\label{soneequiv}
Assume that $\scrB$ is countably thick.
Then $\sone(\seq{\O_n},\scrB)=\sone(\Omega,\scrB)$.
\end{thm}
\begin{proof}
We should verify that the argument in the proof of Theorem \ref{sfinequiv}
works in our case as well.

Let $\seq{\cU_n}$, $\seq{A_n}$, and $\cV_n$ be as in the proof of
Theorem \ref{sfinequiv}.
Apply $\sone(\Omega,\scrB)$ to extract
elements $V_n\in\cV_n$, $n\in\N$, such that $\seq{V_n}\in\scrB$.
For each $n$ and each $m\in A_n$
choose $U_m=V_n$ if $V_n\in\cU_m$, otherwise choose any $U_m\in\cU_m$.
We have enlarged $\seq{V_n}$ by at most countably many open sets.
As $\scrB$ is countably thick, we have that $\seq{U_n}\in\scrB$.
\end{proof}

\section{Examples}

We give some examples for the above results.

\subsection{The Rothberger and Menger properties}
Using our notation, \emph{Rothberger's property $C''$} \cite{ROTH41}
is the property $\sone(\O,\O)$.
In \cite{coc1} it is proved that $\sone(\O,\O)=\sone(\Omega,\O)$.
This implies that $\sone(\Omega,\O)=\sone(\seq{\O_n},\O)$.
Another way to obtain this result is to use Theorem
\ref{soneequiv}, as $\O$ is countably thick.

\emph{Menger's basis property} (introduced in
\cite{MENGER}), was proved by Hurewicz \cite{HURE25}
to be equivalent to the property $\sfin(\O,\O)$.
In \cite{coc1} it is proved that $\sfin(\O,\O)=\sfin(\Omega,\O)$,
so again we have that $\sfin(\Omega,\O)=\sfin(\seq{\O_n},\O)$.

The Rothberger and Menger properties $\sone(\Omega,\O)$ and $\sfin(\Omega,\O)$ are
not provably equivalent, as is witnessed by the canonical Cantor set of reals \cite{coc2}.
Thus, the properties $\sone(\seq{\O_n},\O)$ and $\sfin(\seq{\O_n},\O)$
are not provably equivalent.

\subsection{The Arkhangel'ski\v{i} and Sakai properties}
A space $X$ has the \emph{Arkhangel's\-ki\v{i} property}
\cite{Arkhan} if all finite powers of $X$ have the Menger property
$\sfin(\O,\O)$. In \cite{coc2} it is proved that this is
equivalent to satisfying $\sfin(\Omega,\Omega)$. By Theorem \ref{sfinequiv},
we have that $\sfin(\Omega,\Omega)=\sfin(\seq{\O_n},\Omega)$.

A space $X$ has the \emph{Sakai property} if all finite powers of
$X$ satisfy Rothberger's property $C''$. Sakai \cite{Sakai} proved
that this property is equivalent to $\sone(\Omega,\Omega)$. As
$\Omega$ is countably thick, we have by Theorem \ref{soneequiv}
that $\sone(\Omega,\Omega)=\sone(\seq{\O_n},\Omega)$.

As in the case of Menger and Rothberger, the canonical Cantor set
witnesses that the Arkhangel'ski\v{i} and  Sakai properties
$\sfin(\Omega,\Omega)$ and $\sone(\Omega,\Omega)$ are not provably
equivalent \cite{coc2}. Thus, the properties
$\sone(\seq{\O_n},\Omega)$ and $\sfin(\seq{\O_n},\Omega)$ are not
provably equivalent.

\subsection{The Hurewicz, Gerlits-Nagy $(*)$, and related properties}
$X$ satisfies the \emph{Hurewicz property} (defined in \cite{HURE27})
if for each sequence $\seq{\cU_n}$ of open covers of $X$
there exist finite subsets $\cF_n\sbst\cU_n$, $n\in\N$, such
that $X\sbst\bigcup_n\bigcap_{m>n} \cup\cF_n$
(if $X\nin\seq{\cup\cF_n}$ then this means that $\seq{\cup\cF_n}$
is a $\gamma$-cover of $X$).

To simplify the presentation of the remaining properties, we
introduce the following.
\begin{definition}
~\be
\i
\be
\i A cover $\cU$ of $X$ is \emph{multifinite} if
there exists a partition of $\cU$ into infinitely many
finite covers of $X$.
\i Let $\MF$ denote the collection of all multifinite open covers of $X$.
\ee
\i Fix $\xi\in\{\omega,\gamma,\dots\}$.
\be
\i A cover $\cU$ of $X$ is \emph{$\xi$-groupable}
if it is multifinite, or there exists a partition $\cP$ of $\cU$ into finite sets such that
$\{\cup\cF : \cF\in\cP\}\sm\{X\}$ is a $\xi$-cover of $X$.
\i Let $\O^{\xi\gp}$ denote the collection of all $\xi$-groupable open covers of $X$.
\ee
\ee
\end{definition}

In \cite{coc7} it is proved that
the Hurewicz property
is equivalent to the property $\sfin(\Omega,\O^{\gamma\gp})$.
By Theorem \ref{sfinequiv}, we have that
$X$ has the Hurewicz property if, and only if,
it satisfies $\sfin(\seq{\O_n},\O^{\gamma\gp})$.

In \cite{GN}, Gerlits and Nagy introduced a property called $(*)$.
In \cite{NWS} it is proved that $(*)$ is equivalent to
having the Hurewicz as well as Rothberger properties.
In \cite{coc7} it is proved that
this is equivalent to $\sone(\Omega,\O^{\gamma\gp})$.

\begin{lem}\label{Lgp}
$\O^{\gamma\gp}$ is countably thick.
\end{lem}
\begin{proof}
Assume that $\cU$ is a $\gamma$-groupable cover of $X$,
and let $\cP$ be a partition of $\cU$ witnessing this.
Let $\cV$ be a countable family of open sets.
By shifting to $\cV\sm\cU$ we may assume that
$\cU$ and $\cV$ are disjoint.
As $\cU$ is infinite, $\cP$ is infinite as well;
choose an injection $f:\cV\to\cP$.
Then
$$\tilde\cP = \{f(V)\cup\{V\} : V\in\cV\}\cup(\cP\sm f[\cV])$$
is a partition of $\cU\cup\cV$ witnessing that
this new cover is $\gamma$-groupable.
\end{proof}

\begin{cor}
$\sone(\seq{\O_n},\O^{\gamma\gp})=\sone(\Omega,\O^{\gamma\gp})$.
\end{cor}
Thus a space has the Gerlits-Nagy $(*)$ property
if, and only if, it satisfies $\sone(\seq{\O_n},\O^{\gamma\gp})$.
As the property $(*)$ is not provably equivalent to the
Hurewicz property (this too is witnessed by the Cantor set \cite{coc2},
as $(*)$ implies Rothberger's property \cite{GN}),
we have that $\sone(\seq{\O_n},\O^{\gamma\gp})$ is not provably equivalent to
$\sfin(\seq{\O_n},\O^{\gamma\gp})$.

Now consider the collection $\Omega^{\gpbl}$ of open $\omega$-covers $\cU$ of $X$ such that
there exists a partition $\cP$ of $\cU$ into finite sets such that
for each finite $F\sbst X$ and all but finitely many $\cF\in\cP$,
there exists $U\in\cF$ such that $F\sbst U$.
In \cite{coc7} it is shown that
$X$ satisfies
$\sfin(\Omega,\Omega^{\gpbl})$ if, and only if, all finite powers
of $X$ have the Hurewicz property. By Theorem \ref{sfinequiv},
this property is equivalent to $\sfin(\seq{\O_n},\Omega^{\gpbl})$.
The following observation is what we need to get the
analogous result for the
stronger property $\sone(\Omega,\Omega^{\gpbl})$.
\begin{lem}
$\Omega^{\gpbl}$ is countably thick.
\end{lem}
\begin{proof}
The proof for this is similar to that of Lemma \ref{Lgp}.
\end{proof}
Here too, as all finite powers of the Cantor set $C$ are compact,
we have that $C$ satisfies $\sfin(\Omega,\Omega^{\gpbl})$ but
not $\sone(\Omega,\Omega^{\gpbl})$ (which implies Rothberger's property).
Thus, $\sone(\seq{\O_n},\Omega^{\gpbl})$ is not provably equivalent to
$\sfin(\seq{\O_n},\Omega^{\gpbl})$.

\subsection{A property between Hurewicz and Menger}\label{HureMen}
In \cite{coc1} a property called $\ufin(\Gamma,\Omega)$
is considered, which is intermediate between the
Hurewicz and Menger properties.
This property is a particular case of a general selection
hypothesis. Assume that $\scrA$ and $\scrB$ are collections of covers
of a space $X$. Define the following selection hypothesis \cite{coc1}:
\bi
\item[$\ufin(\scrA,\scrB)$\index{$\ufin(\scrA,\scrB)$}:]
For each sequence $\seq{\cU_n}$ of members of $\scrA$
which do not contain a finite subcover,
there exist finite (possibly empty) subsets $\cF_n\sbst\cU_n$, $n\in\N$,
such that $\seq{\cup\cF_n}\in\scrB$.
\ei
Observe that any countable cover which does not contain a finite
subcover can be turned into a $\gamma$-cover by taking finite unions
\cite{coc2}. Thus for each $\scrA$, $\ufin(\scrA,\scrB)=\ufin(\Gamma,\scrB)$.
The Menger property is equivalent to $\ufin(\Gamma,\O)$,
and the Hurewicz property is equivalent to $\ufin(\Gamma,\Gamma)$.
In \cite{coc2} it is proved that $\ufin(\Gamma,\Omega)$
is not provably equivalent to any of the Hurewicz or Menger properties.

It is proved in \cite{coc8} that
$\ufin(\Gamma,\Omega)$
is equivalent to $\sfin(\Omega,\O^{\omega\gp})$.
By Theorem \ref{sfinequiv},
$X$ satisfies
$\ufin(\Gamma,\Omega)$ if, and only if, it satisfies $\sfin(\seq{\O_n},\O^{\omega\gp})$.

We now treat the stronger property $\sone(\Omega,\O^{\omega\gp})$.
This property was introduced and studied in \cite{coc8}.
In Problem 3 of \cite{coc8} the authors ask whether
this property is strictly stronger than Rothberger's property $\sone(\Omega,\Lambda)$
(this is the same as the usual $\sone(\O,\O)$ \cite{coc1}).
We give a positive answer. It is easy to see (and well known) that
Rothberger's property is closed under taking countable unions.
\begin{thm}
Assuming CH ($\cov(\M)=\c$ is enough),
Rothberger's property does not imply $\sone(\Omega,\O^{\omega\gp})$;
in fact, $\sone(\Omega,\O^{\omega\gp})$ is not even closed
under taking finite unions.
\end{thm}
\begin{proof}
Clearly $\sone(\Omega,\O^{\omega\gp})$ implies
$\sfin(\Omega,\O^{\omega\gp})=\ufin(\Gamma,\Omega)$,
but it is well known that (assuming CH)
Rothberger's property does not imply
$\ufin(\Gamma,\Omega)$ \cite{coc2}.

Moreover, in \cite{huremen2} it is shown
that CH (or even just $\cov(\M)=\c$) implies
that no property between $\sone(\Omega,\Omega)$ and
$\ufin(\Gamma,\Omega)$ (inclusive) is closed under taking finite
unions. But $\sone(\Omega,\O^{\omega\gp})$ lies
between these properties.
\end{proof}

As in the case of $\sone(\Omega,\O^{\gamma\gp})$ which is equivalent
to $\ufin(\Gamma,\Gamma)\cap \sone(\O,\O)$, we have that
the new property $\sone(\Omega,\O^{\omega\gp})$ can also be characterized
in terms of the more classical properties.
\begin{thm}
$\sone(\Omega,\O^{\omega\gp}) = \ufin(\Gamma,\Omega)\cap\sone(\O,\O)$.
\end{thm}
\begin{proof}
We have seen that $\sone(\Omega,\O^{\omega\gp})$ implies $\ufin(\Gamma,\Omega)$
and Rothberger's property $\sone(\O,\O)$. To prove the other implication,
we use the result of \cite{coc8}, that $\ufin(\Gamma,\Omega)$ implies
$\Lambda=\O^{\omega\gp}$. As $\sone(\O,\O)=\sone(\Omega,\Lambda)$,
$\Lambda=\O^{\omega\gp}$ and $\sone(\O,\O)$ imply $\sone(\Omega,\O^{\omega\gp})$.
\end{proof}

\begin{lem}
$\O^{\omega\gp}$ is countably thick.
\end{lem}
\begin{proof}
Assume that $\cU$ is an $\omega$-groupable cover of $X$,
and let $\cP$ be a partition of $\cU$ witnessing this.
Let $\cV$ be a countable family of open sets.
We may assume that $\cU$ and $\cV$ are disjoint.
Let $\tilde\cP$ be any partition of $\cV$ into
finite sets.
Then $\cP\cup\tilde\cP$
is a partition of $\cU\cup\cV$ witnessing that
this new cover is $\omega$-groupable.
\end{proof}

\begin{cor}
$\sone(\seq{\O_n},\O^{\omega\gp})=\sone(\Omega,\O^{\omega\gp})$.
\end{cor}
Here again, Cantor's set witnesses that the properties $\sone(\seq{\O_n},\O^{\omega\gp})$ and
$\sfin(\seq{\O_n},\O^{\omega\gp})$ are not provably equivalent.

\section{$\tau$-covers}
An open cover $\cU$ of $X$ is a \emph{$\tau$-cover} of $X$
if it is a large cover, and for each $x,y\in X$,
one of the sets $\{U\in\cU : x\in U\mbox{ and }y\nin U\}$
or $\{U\in\cU : y\in U\mbox{ and }x\nin U\}$ is finite.
The notion of $\tau$-covers was introduced in \cite{tau},
and incorporated into the framework of selection principles in
\cite{tautau}.

Let $\Tau$ denote the collection of open $\tau$-covers of $X$.
Then $\Gamma\sbst\Tau\sbst\Omega$, therefore
$\sone(\Omega,\Gamma)$ implies $\sone(\Omega,\Tau)$, which implies
$\sfin(\Omega,\Tau)$. It is not known whether any two of
these properties are equivalent.

By Theorem \ref{sfinequiv}, we have that
$\sfin(\seq{\O_n},\Tau)=\sfin(\Omega,\Tau)$.
We have only a guess for the situation in the
remaining case.
\begin{conj}\label{conj}
It is consistent that $\sone(\seq{\O_n},\Tau)\neq\sone(\Omega,\Tau)$.
\end{conj}
Observe that, as $\sone(\Omega,\Tau)$ implies Rothberger's property
$\sone(\O,\O)$, we have by the consistency of Borel's conjecture that
the word ``consistent'' cannot be replaced by ``provable'' in
Conjecture \ref{conj}.

$\tau^*$-covers are a variation of $\tau$-covers which is
easier to work with.
For a cover $\cU=\seq{U_n}$ of $X$ and an element $x\in X$, write
$$x_\cU = \{n : x\in U_n\}.$$
According to \cite{tautau},
A cover $\cU$ of $X$ is a \emph{$\tau^*$-cover} of $X$ if it is
large, and for each $x\in X$ there exists an infinite subset $\hat x_\cU$ of $x_\cU$
such that the sets $\hat x_\cU$, $x\in X$, are linearly quasiordered by
$\as$ ($A\as B$ means that $A\sm B$ is finite).
If $\cU$ is a countable $\tau$-cover, then
by setting $\hat x_\cU = x_\cU$ for each $x\in X$ we see that
it is a $\tau^*$-cover. The converse is not
necessarily true.
Let $\Tau^*$ denote the collection of all countable open
$\tau^*$-covers of $X$. Then $\Tau\sbst\Tau^*\sbst\Omega$.

\begin{lem}
$\Tau^*$ is countably and finitely thick.
\end{lem}
\begin{proof}
Assume that $\cU=\seq{U_n}\in \Tau^*$, and
let $\hat x_\cU$, $x\in X$, be witnesses for that.
Let $\cV$ be a countable family of open sets.
Assume that $\cV$ is infinite and disjoint
from $\cU$, and let $\seq{V_n}$ be a bijective enumeration
of $\cV$.
Enumerate $\cU\cup\cV$ by
$\seq{W_n}$ where $W_n = U_n$ if $n$ is even
and $W_n = V_n$ otherwise.
Then the subsets $2\hat x_\cU$ of $x_{\cU\cup\cV}$, $x\in X$, witness that
$\cU\cup\cV\in\Tau^*$.
The case that $\cV$ has a finite cardinality $k$ is treated similarly.

To see that $\Tau^*$ is finitely thick it remains to
verify the first requirement in the definition of finitely thick covers.
In \cite{tautau} we prove something stronger: If $\cU\in\Tau^*$ refines
a countable cover $\cV$, then $\cV\in\Tau^*$.
\end{proof}

\begin{cor}
$\sone(\seq{\O_n},\Tau^*)=\sone(\Omega,\Tau^*)$.
\end{cor}
The last corollary can be contrasted with Conjecture \ref{conj}.

\section{Iliadis' question}

In the
\emph{Lecce Workshop on Coverings, Selections and Games in Topology}
(June 2002), Stavros Iliadis asked whether we get new properties
if we consider the generalized selection principles of the
form $\sone(\seq{\scrA_n},\scrB)$ and $\sfin(\seq{\scrA_n},\scrB)$.
We check the cases where the first coordinate
is any sequence of elements from the set
$$\cC=\{\O,\Lambda,\Omega,\Tau^*,\Tau,\Gamma\}\cup\{\O_n : n\in\N\}.$$

\begin{lem}\label{gammasubseq}
For any increasing sequence $\seq{k_n}$,
$\sone(\seq{\scrA_n},\Gamma)$ implies $\sone(\seq{\scrA_{k_n}},\Gamma)$,
and $\sfin(\seq{\scrA_n},\Gamma)$ implies $\sfin(\seq{\scrA_{k_n}},\Gamma)$.
\end{lem}
\begin{proof}
Assume that $\cU_{k_n}\in\scrA_{k_n}$.
For each $m\nin\seq{k_n}$ use Lemma \ref{alwaysgamma} in Appendix \ref{toostrong}
to choose an element $\cU_m\in\Gamma$.

Apply $\sone(\seq{\scrA_n},\Gamma)$
to the sequence
$\seq{\cU_n}$ to obtain elements $U_n\in\cU_n$
such that $\seq{U_n}$ is a
$\gamma$-cover of $X$.
Then $\seq{U_{k_n}}$ is a
$\gamma$-cover of $X$,
and for each $n$, $U_{k_n}\in\cU_{k_n}$.

The proof for $\sfin$ is similar.
\end{proof}

\begin{cor}\label{s1gamma}
Assume that $\seq{\scrA_n}$ is a sequence of elements of $\cC$.
Then:
\be
\i If some $\scrA\in\{\O, \Lambda\}\cup\{\O_n : n\in\N\}$
occurs infinitely often in the sequence $\seq{\scrA_n}$,
then $\sone(\seq{\scrA_n},\Gamma)$ implies $\sone(\scrA,\Gamma)$, which
is false for a nontrivial $X$.
\i If $\scrA_n\in\{\O, \Lambda\}$ for only finitely many $n$
and there exists an increasing sequence $k_n$ such that $\seq{\O_{k_n}}$ is a
subsequence of $\seq{\scrA_n}$,
then $\sone(\seq{\scrA_n},\Gamma)=\sone(\seq{\O_n},\Gamma)$ (strong $\gamma$-set).
\i  If $\scrA_n\in\{\O, \Lambda\}\cup\{\O_n : n\in\N\}$ for only finitely many $n$
and $\cU_n=\Omega$ for infinitely many $n$, then
$\sone(\seq{\scrA_n},\Gamma)=\sone(\Omega,\Gamma)$.
\i If $\scrA_n\in\{\O, \Lambda,\Omega\}\cup\{\O_n : n\in\N\}$ for only finitely many $n$
and $\cU_n=\Tau^*$ for infinitely many $n$, then
$\sone(\seq{\scrA_n},\Gamma)=\sone(\Tau^*,\Gamma)$.
\i If $\scrA_n\in\{\O, \Lambda,\Omega,\Tau^*\}\cup\{\O_n : n\in\N\}$ for only finitely many $n$
and $\cU_n=\Tau$ for infinitely many $n$, then
$\sone(\seq{\scrA_n},\Gamma)=\sone(\Tau,\Gamma)$.
\i If $\scrA_n\in\{\O, \Lambda,\Omega,\Tau^*,\Tau\}\cup\{\O_n : n\in\N\}$
for only finitely many $n$, then
$\sone(\seq{\scrA_n},\Gamma)=\sone(\Gamma,\Gamma)$.
\ee
The analogous assertions for $\sfin$ also hold.
\end{cor}
\begin{proof}
We will use Lemma \ref{gammasubseq}.

(1) follows, using the results of
Appendix \ref{toostrong}.

To prove (2), observe that in this case, each $\cU_n$
is a subset of $\O_{k_n}$ for some $k_n$, such that
$k_n\to\infty$, so that
Theorem \ref{sgnicedef} applies.

(3), (4), (5), and (6) follow from Lemma \ref{gammasubseq}.
\end{proof}

\begin{lem}\label{supU}
\be
\i If $\scrA_n\spst\scrA$ for all but finitely many $n$,
and $\scrB$ is closed under removing a finite subset, then
$\sone(\seq{\scrA_n},\scrB)$ implies $\sone(\scrA,\scrB)$.
\i If $\scrA$ occurs infinitely often in the sequence
$\seq{\scrA_n}$, and $\scrB$ is countably thick, then
$\sone(\scrA,\scrB)$ implies $\sone(\seq{\scrA_n},\scrB)$.
\i The same assertions hold for $\sfin$ (where in (2)
countable thickness is not needed).
\ee
\end{lem}
\begin{proof}
(1) Assume that $X$ satisfies $\sone(\seq{\scrA_n},\scrB)$.
We will show that  $X$ satisfies $\sone(\scrA,\scrB)$.
Fix $m$ such that that for all $n\geq m$, $\scrA_n\spst\scrA$.
Assume that $\seq{\cU_n}$ is such that $\cU_n\in\scrA$ for all
$n$.
By Lemma \ref{alwaysgamma}, there exists an open $\gamma$-cover
$\cV$ of $X$. Define a sequence $\seq{\cV_n}$ by
$\cV_n=\cV$ for $n< m$ and $\cV_n=\cU_{n-m}$ otherwise.
By $\sone(\seq{\scrA_n},\scrB)$, there exist elements
$V_n\in\cV_n$ such that $\seq{V_n}\in\scrB$.
As $\scrB$ is closed under removing a finite subset,
$\{V_n\}_{n\ge m}\in\scrB$ and for each $n\ge m$, $V_n\in\cU_{n-m}$.

(2) Let $\seq{k_n}$ be an increasing enumeration of $\{n : \scrA_n=\scrA\}$,
and let $\seq{\cU_n}$ be such that $\cU_n\in\scrA_n$ for all $n$.
Apply $\sone(\scrA,\scrB)$
to $\cU_{k_n}$ to obtain elements $U_{k_n}\in\cU_{k_n}$
such that $\seq{U_{k_n}}$ is a member of $\scrB$.
From the remaining covers $\cU_n$ choose any element $U_n$.
As $\scrB$ is countably thick, $\seq{U_{k_n}}$ is a member of $\scrB$
as well.

(3) is similar.
\end{proof}

The collections $\Lambda$, $\Omega$, $\Tau^*$ and $\Tau$ are all
countably thick and closed under removing a finite subset.
Thus, if $\scrB$ is any of these, then we get
$\sone(\seq{\scrA_n},\scrB)=\sone(\scrA,\scrB)$ in Lemma \ref{supU}.

\begin{cor}\label{soneallcases}
Assume that $\seq{\scrA_n}$ is a sequence of elements of $\cC$,
and $\scrB\in\{\Lambda, \Omega\}$.
Then:
\be
\i If there exist infinitely many $n$ such that $\scrA_n=\Gamma$,
then $\sone(\seq{\scrA_n},\scrB)=\sone(\Gamma,\scrB)$.
\i If there exist only finitely many $n$ such that $\scrA_n=\Gamma$,
and there exist infinitely many $n$ such that $\scrA_n=\Tau$,
then $\sone(\seq{\scrA_n},\scrB)=\sone(\Tau,\scrB)$.
\i If there exist only finitely many $n$ such that $\scrA_n\in\{\Tau,\Gamma\}$,
and there exist infinitely many $n$ such that $\scrA_n=\Tau^*$,
then $\sone(\seq{\scrA_n},\scrB)=\sone(\Tau^*,\scrB)$.
\i If there exist only finitely many $n$ such that $\scrA_n\in\{\Tau^*,\Tau,\Gamma\}$,
and there exist infinitely many $n$ such that $\scrA_n=\Omega$,
then $\sone(\seq{\scrA_n},\scrB)=\sone(\Omega,\scrB)$.
\i If there exist only finitely many $n$ such that $\scrA_n\in\{\Omega,\Tau^*,\Tau,\Gamma\}$,
and there exists an increasing sequence $\seq{k_n}$ such that $\scrA_{k_n}\sbst\O_n$ for all $n$,
then $\sone(\seq{\scrA_n},\scrB)=\sone(\Omega,\scrB)$.
\i If there exists no increasing sequence $\seq{k_n}$ such that $\scrA_{k_n}\sbst\O_n$ for all $n$,
and $\Lambda$ occurs infinitely often in $\seq{\scrA_n}$, then
$\sone(\seq{\scrA_n},\scrB)=\sone(\Lambda,\scrB)$ (which is Rothberger's property if
$\scrB\in\{\O,\Lambda\}$ and trivial otherwise).
\i If for some $m$ $\scrA_n\spst\O_m$ for almost all $n$, then
$\sone(\seq{\scrA_n},\Lambda)$ is trivial.
\ee
The analogous assertions for $\sfin$ also hold.
\end{cor}

These results and related arguments should show that no new property is
introduced by the generalized selection principles $\sone(\seq{\scrA_n},\scrB)$ and
$\sfin(\seq{\scrA_n},\scrB)$, except for the strong $\gamma$-property
$\sone(\seq{\O_n},\Gamma)$ and, perhaps, $\sone(\seq{\O_n},\Tau)$.

\bigpart{Part II: Game theory}
In this section we give new game theoretic characterizations to
most of the properties considered in the previous sections.
Although these characterizations are suggested by the results
of the earlier sections, their proofs are not as trivial.

\section{Selection games and strategies}
Each selection principle has a naturally associated game.
In the game $\gone(\scrA,\scrB)$ ONE chooses in the $n$th inning an element
$\cU_n$ of $\scrA$ and then TWO responds by choosing $U_n\in \cU_n$.
They play an inning per natural number.
A play $(\cU_0, U_0, \cU_1, U_1,\dots)$ is won by TWO if
$\seq{U_n}\in\scrB$; otherwise ONE wins.
The game $\gfin(\scrA,\scrB)$ is played similarly, where TWO responds with
finite subsets $\cF_n\sbst\cU_n$ and wins if $\Union_{n\in\N}\cF_n\in\scrB$.

Observe that if ONE does not have a winning strategy in
$\gone(\scrA,\scrB)$ (respectively, $\gfin(\scrA,\scrB)$), then $\sone(\scrA,\scrB)$
(respectively, $\sfin(\scrA,\scrB)$) holds.
The converse is not always true; when it is true,
the game is a powerful tool for studying the combinatorial
properties of $\scrA$ and $\scrB$ -- see, e.g., \cite{coc7, coc8}, and references therein.

It is therefore appealing to try and study the generalized games
associated with $\sone(\seq{\scrA_n},\scrB)$ and $\sfin(\seq{\scrA_n},\scrB)$.
\begin{definition}
Define the following games between two players, ONE and TWO,
which have an inning per natural number.
\bi
\i[$\gone(\seq{\scrA_n},\scrB)$:]
In the $n$th inning, ONE chooses an element $\cU_n\in\scrA_n$, and
TWO responds with an element $U_n\in\cU_n$. TWO wins if
$\seq{U_n}\in\scrB$; otherwise ONE wins.
\i[$\gfin(\seq{\scrA_n},\scrB)$:]
In the $n$th inning, ONE chooses an element $\cU_n\in\scrA_n$, and
TWO responds with a finite subset $\cF_n$ of $\cU_n$. TWO wins if
$\Union_{n\in\N}\cF_n\in\scrB$; otherwise ONE wins.
\ei
\end{definition}

Some terminological conventions will be needed to simplify the
proofs of the upcoming results.
A \emph{strategy} $F$ for ONE in a game $\gone(\seq{\scrA_n},\scrB)$
can be identified with a tree
of covers in the following way.
Let $\cU_{\<\>}:=F(X)$ be the first move of ONE.
Enumerate the elements of $\cU_{\<\>}$ as $\seq{U_{\<n\>}}$.
Having defined $\cU_{\sigma}=\seq{U_{\sigma\cat \<n\>}}$,
define for each $m$
$$\cU_{\sigma\cat \<m\>}:= F(\cU_{\<\>}, U_{\sigma\| 1},\dots,\cU_\sigma,U_{\sigma\cat\<m\>}),$$
and fix an enumeration $\seq{U_{\sigma\cat\<m,n\>}}$ of $\cU_{\sigma\cat \<m\>}$.
Let $\N^*$ denote the collection of all finite sequences of natural numbers.

Similarly, a strategy $F$ for ONE in a game $\gfin(\seq{\scrA_n},\scrB)$
can be identified with a tree covers where the sequences $\sigma$ are of
\emph{finite sets} of
natural numbers rather than natural numbers.
Let $[\N]^*$ denote the collection of all finite sequences of finite
sets of natural numbers.

We will say that a collection of covers $\scrA$ is
\emph{dense} in a strategy $F$ for ONE in a game of type $\gone$
if for each $\sigma\in\N^*$
there exists $\eta\in\N^*$ which extends
$\sigma$, and such that $\cU_\eta\in\scrA$, that is, $\{\eta\in\N^* : \cU_\eta\in\scrA\}$ is
dense in $\N^*$.
Similarly, we say that $\scrA$ is dense in a strategy
$F$ for ONE in a game of type $\gfin$ if $\{\eta\in[\N]^* : \cU_\eta\in\scrA\}$
is dense in $[\N]^*$.

\begin{lem}[Density lemma]\label{densecover}
Assume that $\scrB$ is countably thick,
ONE has a winning strategy $F$
in $\gone(\seq{\scrA_n},\scrB)$,
and $\scrA$ is dense in $F$.
Then ONE has a winning strategy in $\gone(\scrA,\scrB)$.
The analogous assertion for $\gfin$ also holds.
\end{lem}
\begin{proof}
Fix some well-ordering on the collection $\N^*$ of
all finite sequences of natural numbers, and
let $\{\cU_\sigma\}_{\sigma\in\N^*}$ be the
tree of covers associated with $F$.
Define a function $\pi:\N^*\to\N^*$
as follows:
\be
\i Let $\pi(\<\>)$ be the first member of $\N^*$
such that $\cU_{\pi(\<\>)}\in\scrA$.
\i For each $n$ let $\pi(\<n\>)\in\N^*$
be the first extension of $\pi(\<\>)\cat\<n\>$
such that $\cU_{\pi(\<n\>)}\in\scrA$.
\i In general, for each $\sigma\in\N^*$ and each $n$
let $\pi(\sigma\cat\<n\>)$ be the first extension
of $\pi(\sigma)\cat\<n\>$ such that $\cU_{\pi(\sigma\cat\<n\>)}\in\scrB$.
\ee
For each $\sigma\in\N^*$ define $\tilde\cU_\sigma=\cU_{\pi(\sigma)}$, and
set $\tilde U_{\sigma\cat\<n\>}=U_{\pi(\sigma)\cat\<n\>}$ for each $n$.
Let $\tilde F$ be the strategy associated with
$\{\tilde\cU_\sigma\}_{\sigma\in\N^*}$.
Then $\tilde F$ is a strategy for ONE in $\gone(\scrA,\scrB)$.

We claim that $\tilde F$ is a winning strategy for
ONE in $\gone(\scrA,\scrB)$.
Assume otherwise, and let $f\in\NN$ be such that the play
$$(\tilde\cU_{\<\>}, \tilde U_{f\| 1}, \tilde\cU_{f\|1}, \tilde U_{f\| 2},\dots)$$
against the strategy $\tilde F$ is lost by ONE, that is, $\seq{U_{f\| n}}\in\scrB$.
Define
$$\sigma_0=\pi(\<\>), \sigma_1=\pi(\<f(0)\>), \dots, \sigma_{n+1}=\pi(\sigma_n\cat \<f(n)\>),\dots$$
and take $g=\Union_{n\in\N}\sigma_n$.
Then
$$(\cU_{\<\>}, U_{g\| 1}, \cU_{g\|1}, U_{g\| 2},\dots)$$
is a play in the game $\gone(\seq{\scrA_n},\scrB)$ according to the strategy $F$,
and $\seq{U_{f\| n}}$ is a subsequence of $\seq{U_{g\| n}}$.
As $\scrB$ is countably thick, we have that $\seq{U_{g\| n}}\in\scrB$ as well,
so this game is lost by ONE, a contradiction.

The proof for $\gfin$ is similar.
\end{proof}

\begin{rem}
For each $\sigma\in\N^*$, we can modify the definition of $\pi$ in the
proof of Lemma \ref{densecover} so that $\pi(\<\>)$ extends $\sigma$.
Consequently, it is enough to assume that $\{\eta : \cU_\eta\in\scrA\}$
is dense below $\sigma$ (with respect to the order of reverse inclusion)
for \emph{some} $\sigma\in\N^*$.
In other words, if ONE does not have a winning strategy in
$\gone(\scrA,\scrB)$ but has a winning strategy in $\gone(\seq{\scrA_n},\scrB)$,
then $\{\eta : \cU_\eta\in\scrA\}$ is \emph{nowhere dense} in $\N^*$.
\end{rem}

\section{Reductions among $\gfin$ strategies}

Following is a surprising result.
It implies that if $\MF\sbst\scrB$ and
ONE could win the game $\gfin(\seq{\O_n},\scrB)$,
then he could win $\gfin(\Lambda,\scrB)$ as well.

\begin{thm}[$\Lambda$-less strategies]\label{denseLam}
Assume that $F$ is a strategy for ONE in a game $\gfin(\seq{\O_n},\scrB)$,
and $\Lambda$ is not dense in $F$.
Then $F$ is not even a winning strategy for ONE in the game
$\gfin(\seq{\O_n},\MF)$.
\end{thm}
\begin{proof}
Assume that $F$ is a winning strategy for ONE in the game $\gfin(\seq{\O_n},\MF)$.
Let $\{\cU_\sigma\}_{\sigma\in[\N]^*}$ be the covers tree associated
with $F$, and
choose $\sigma\in[\N]^*$ such that for all $\eta$ extending $\sigma$,
$\tilde\cU_\eta$ is not large.
Modify $F$ so that its first move is $\tilde\cU_\eta$
(that is, the strategy determined by the subtree $\{\sigma : \eta\sbst\sigma\}$
of $[\N]^*$).
This is still a winning strategy for ONE (otherwise
TWO can begin with a sequence of moves which will
force ONE into $\tilde\cU_\eta$ and then defeat him).
We may therefore assume that no element in the strategy $F$ is large.

\begin{lem}
Every $n+1$ cover of a space $X$ which is not large contains a finite
$n$-cover of $X$.
\end{lem}
\begin{proof}
Assume that $\cU$ is an $n+1$-cover of $X$ which is not large.
Then there exists
$x\in X$ such that the set $\cF=\{U\in\cU : x\in U\}$
is finite.
Now, as $\cU$ is an $n+1$-cover of $X$,
for each $n$-element subset $F$ of $X$
there exists $U\in\cU$ such that $F\cup\{x\}\sbst U$, and
therefore $U\in\cF$ and $F\sbst U$.
\end{proof}
We may therefore modify the strategy $F$ (by thinning out its covers)
so that all covers in this strategy are finite.
As this only restricts the possible moves of TWO, this is
still a winning strategy for ONE in the game $\gfin(\seq{\O_n},\MF)$.

In particular, no cover in the strategy $F$ is an $\omega$-cover of $X$.
\begin{lem}\label{disjointlemma}
Assume that $\seq{\cU_n}$ is a sequence of $n$-covers of
$X$ which are not $\omega$-covers of $X$.
Then there exists an increasing sequence $\seq{k_n}$
and pairwise disjoint subsets $\tilde\cU_{k_n}$ of $\cU_{k_n}$ such
that each $\tilde\cU_{k_n}$ is an $n$-cover of $X$.
\end{lem}
\begin{proof}
For each $n$ let $F_n$ be a finite subset of
$X$ witnessing that $\cU_n$ is not an $\omega$-cover of $X$.
Observe that if $\cU$ is a $k+l$-cover of $X$ and
$F$ is a $k$-element subset of $X$. Then
$\{U\in\cU : F\sbst U\}$ is an $l$-cover of $X$.

Let $\tilde\cU_1=\cU_1$. Set $k_1=|F_1|+1$.
Then $\tilde\cU_2=\{U\in\cU_{k_1} : F_1\sbst U\}$ is a cover
of $X$ disjoint from $\tilde\cU_1$.
Assume that we have defined $\tilde\cU_{k_1},\dots,\tilde\cU_{k_{n-1}}$.
Let $k_n=|\Union_{i<n}F_{k_i}|+n$, and
choose
$$\tilde\cU_{k_n} = \{U\in\cU_{k_n} : \Union_{i<k}F_{k_i}\sbst U\}.$$
Then $\tilde\cU_{k_n}$ is an $n$-cover of $X$,
$\tilde\cU_{k_n}\sbst\cU_{k_n}$, and $\tilde\cU_{k_n}\cap\tilde\cU_{k_i}=\emptyset$
for all $i<n$.
\end{proof}
Thus, by the methods of Lemma \ref{densecover},
we may refine the strategy $F$ so that all
its covers are (finite and) disjoint.
Again, as the new strategy restricts the moves of
TWO, it is still a winning strategy in the game
$\gfin(\seq{\O_n},\MF)$.
But in this situation TWO can choose the whole
cover in each inning, making its confident way to
a victory in the game $\gfin(\seq{\O_n},\MF)$,
a contradiction.
\end{proof}

We now give some applications of Theorem \ref{denseLam}.

\begin{thm}\label{MenGame}
For a space $X$, the following are equivalent.
\be
\i $X$ has the Menger property,
\i ONE does not have a winning strategy in $\gfin(\Omega,\Lambda)$,
\i ONE does not have a winning strategy in $\gfin(\Lambda,\Lambda)$; and
\i ONE does not have a winning strategy in $\gfin(\seq{\O_n},\Lambda)$.
\ee
\end{thm}
\begin{proof}
Hurewicz \cite{HURE25} proved that the Menger property is
equivalent to ONE not having a winning strategy in $\gfin(\O,\O)$.
Using this and the method in Theorem 3 of \cite{OpPar},
one shows that $1\Iff 3$.
Now, $3\Impl 2$, and $2$ implies Menger's property $\sfin(\Omega,\Lambda)$.
Thus $1\Iff 2\Iff 3$. Clearly $4\Impl 2$.
To see that $3\Impl 4$, assume that $F$ is a strategy for ONE in
$\gfin(\seq{\O_n},\Lambda)$. By $3$ and  Lemma \ref{densecover},
$\Lambda$ is not dense in $F$. By Theorem \ref{denseLam},
$F$ is not a winning strategy for ONE in  $\gfin(\seq{\O_n},\Lambda)$.
\end{proof}

\begin{thm}\label{HureGame}
For a space $X$, the following are equivalent.
\be
\i $X$ has the Hurewicz property,
\i ONE does not have a winning strategy in $\gfin(\Omega,\O^{\gamma\gp})$,
\i ONE does not have a winning strategy in $\gfin(\Lambda,\O^{\gamma\gp})$; and
\i ONE does not have a winning strategy in $\gfin(\seq{\O_n},\O^{\gamma\gp})$.
\ee
\end{thm}
\begin{proof}
The equivalence $1\Iff 2$ is established in Theorem 12 of \cite{coc7}.
It is clear that $3\Impl 2$ and $4\Impl 2$.

$1\Impl 3$:
This is proved like the proof of $1\Impl 2$
(see Theorem 12 of \cite{coc7}),
as $\Lambda$ is closed under removing a finite subset.

$3\Impl 4$:
Assume that ONE has a winning strategy $F$
in the game $\gfin(\seq{\O_n},\O^{\gamma\gp})$.
Then by Theorem \ref{denseLam},
$\Lambda$ is dense in $F$, and by Lemma \ref{densecover}, ONE has a winning strategy in
the game $\gfin(\Lambda,\O^{\gamma\gp})$.
\end{proof}

The last $\gfin$ game we consider is the one associated to the property
$\ufin(\Gamma,\Omega)$ from Subsection \ref{HureMen}.

\begin{thm}\label{HureMenGame}
For a space $X$, the following are equivalent.
\be
\i $X$ satisfies $\ufin(\Gamma,\Omega)$,
\i ONE does not have a winning strategy in $\gfin(\Omega,\O^{\omega\gp})$,
\i ONE does not have a winning strategy in $\gfin(\Lambda,\O^{\omega\gp})$; and
\i ONE does not have a winning strategy in $\gfin(\seq{\O_n},\O^{\omega\gp})$.
\ee
\end{thm}
\begin{proof}
$1\Iff 2$ is proved in Theorem 13 of \cite{coc8}.

$2\Impl 3$: Assume 2. Then ONE does not have
a winning strategy in the game $\gfin(\Omega,\Lambda)$.
By Theorem \ref{MenGame}, ONE does not have
a winning strategy in the game $\gfin(\Lambda,\Lambda)$.
According to Lemma 11 of \cite{coc8}, 1 (which is implied by 2) implies
that each large cover of $X$ is $\omega$-groupable,
that is, $\Lambda=\O^{\omega\gp}$ for $X$.
Thus ONE does not have
a winning strategy in the game $\gfin(\Lambda,\O^{\omega\gp})$.

$2\Impl 4$ is proved similarly.
\end{proof}

The following problem remains open.
\begin{prob}
Is the Arkhangel'ski\v{i} Property
$\sfin(\Omega,\Omega)$ equivalent to
ONE not having a winning strategy in $\gfin(\{\O_n\},\Omega)$ ?
\end{prob}

\section{Reductions among $\gone$ strategies}

We now turn to $\gone$-games.
To deal with these, we need some more terminology and tools.
Assume that $F$ is a strategy for ONE in a $\gone$-game.
The \emph{$\omega$-strategy} $F_\omega$ associated to
$F$ is the strategy defined as follows.
Let $\{\cU_\sigma\}_{\sigma\in\N^*}$ be the covers tree
associated to $F$.
Fix a bijection $\phi:\N\to\N^*$.
For each $\sigma\in\N^*$ of length $k$ let
$$\phi(\sigma)=\phi(\sigma(0))\cat\phi(\sigma(1))\cat\dots\cat\phi(\sigma(k-1)).$$
For each $n$
define $\hat U_{\<n\>}=\cup\{U_{\eta\|1},U_{\eta\|2}\dots,U_{\eta}\}$ where $\eta=\phi(n)$,
and set $\hat\cU_{\<\>} = \seq{\hat U_{\<n\>}}$. In general,
for each $\sigma\in\N^*$ and each $n$ let $\eta=\phi(n)$, and define
$$\hat U_{\sigma\cat\<n\>}=\cup\{U_{\phi(\sigma)\cat\eta\|1},U_{\phi(\sigma)\cat\eta\|2}\dots,
U_{\phi(\sigma)\cat\eta}\}.$$
Set $\hat\cU_\sigma=\seq{\hat U_{\sigma\cat\<n\>}}$.

As we have required that $X$ is not a member of any cover we consider,
$F_\omega$ need not be a strategy for ONE.
We will say that $X$ is \emph{$\omega$-dense} in $F$
if the set $\{\sigma : \hat U_\sigma=X\}$ is dense in $\N^*$.
\begin{lem}\label{w-dense}
Assume that for each $\sigma$,
$\cU_\sigma$ is disjoint from its past $\{U_{\sigma\|1},U_{\sigma\|2},\dots,U_{\sigma}\}$.
If $X$ is $\omega$-dense in $F$, then there exists a game according to this strategy
where the moves of TWO constitute a groupable large cover of $X$.
\end{lem}
\begin{proof}
In the covers tree of $F$ there exists a path with infinitely many disjoint
intervals which constitute a finite cover of $X$.
\end{proof}

\begin{lem}\label{not-w-dense}
Fix $\scrB\in\{\Lambda,\O^{\omega\gp},\O^{\gamma\gp}\}$.
Assume that $F$ is a strategy for ONE in $\gone(\seq{\scrA_n},\scrB)$
such that $X$ is not $\omega$-dense in $F$,
and for each $\sigma$,
$\cU_\sigma$ is disjoint from its past $\{U_{\sigma\|1},\allowbreak U_{\sigma\|2},\allowbreak \dots,
\allowbreak U_{\sigma}\}$.
If $F_\omega$ is not a winning strategy for
ONE in the game $\gone(\Omega,\scrB)$,
then $F$ is not a winning strategy for ONE in the game $\gone(\seq{\scrA_n},\scrB)$.
\end{lem}
\begin{proof}
Any move of TWO in $F_\omega$ can be translated to a finite sequence of moves
for TWO in $F$, replacing each $\hat U_{\sigma\cat n}$ chosen by TWO with the
elements
$U_{\phi(\sigma)\cat\eta\|1},\allowbreak U_{\phi(\sigma)\cat\eta\|2},\dots,\allowbreak
U_{\phi(\sigma)\cat\eta}$ where $\eta=\phi(n)$.
It is easy to see, by disjointness from the past, that this disassemblying
preserves being a member of $\scrB$ for $\scrB\in\{\Lambda,\O^{\omega\gp},\allowbreak \O^{\gamma\gp}\}$.
\end{proof}

For shortness, we give the characterizations for
the Rothberger,  Gerlits-Nagy $(*)$, and $\sone(\Omega,\O^{\omega\gp})$ properties simultaneously.
\begin{thm}\label{GoneGame}
Fix $\scrB\in\{\Lambda,\O^{\omega\gp},\O^{\gamma\gp}\}$.
For a space $X$, the following are equivalent.
\be
\i $X$ satisfies $\sone(\Omega,\scrB)$,
\i ONE does not have a winning strategy in $\gone(\Omega,\scrB)$,
\i ONE does not have a winning strategy in $\gone(\Lambda,\scrB)$; and
\i ONE does not have a winning strategy in $\gone(\seq{\O_n},\scrB)$.
\ee
\end{thm}
\begin{proof}
In Theorem 3 of \cite{OpPar} it is proved that $1\Iff 3$ for $\scrB=\Lambda$.
In Theorem 12 of \cite{coc7} it is proved that $1\Iff 2$ for $\scrB=\O^{\gamma\gp}$,
and in Theorem 15 of \cite{coc8} this is proved for $\scrB=\O^{\omega\gp}$.

$2\Impl 3$:
Assume that $F$ is a winning strategy for ONE in
$\gone(\Lambda,\scrB)$.
Modify the covers tree by removing from each
$\cU_\sigma$ its past $\{U_{\sigma\|1},U_{\sigma\|2},\dots,U_{\sigma}\}$.
Then $F$ is still a winning strategy for ONE.
By Lemma \ref{w-dense}, $X$ is not $\omega$-dense in $F$,
and by Lemma \ref{not-w-dense} we get that $F_\omega$ is a winning strategy
for ONE in $\gone(\Omega,\scrB)$.

$2\Impl 4$: Assume that $F$ is a winning strategy for ONE in
$\gone(\seq{\O_n},\scrB)$.
If $\Omega$ is dense in $F$ then by Lemma \ref{densecover}
ONE has a winning strategy in $\gone(\Omega,\scrB)$.
Otherwise, by Lemma \ref{disjointlemma} we may assume that
the covers in each branch of the strategy $F$ are disjoint.
By Lemma \ref{w-dense}, $X$ is not $\omega$-dense in $F$,
and by Lemma \ref{not-w-dense} we get that $F_\omega$ is a winning strategy
for ONE in $\gone(\Omega,\scrB)$.
\end{proof}

\begin{rem}
The characterizations of Rothberger's property
using $\O$ instead of $\Lambda$
are much more simple to deal with:
Pawlikowski \cite{PAW} proved that Rothberger's property
$\sone(\O,\O)$ is equivalent to ONE not having a winning
strategy in $\gone(\O,\O)$.
As $\sone(\O,\O)=\sone(\Omega,\O)$, we get that
Rothberger's property is equivalent to ONE not having a winning
strategy in  $\gone(\seq{\O_n},\O)$.
\end{rem}

We now treat the remaining $\gone$-games: $\gone(\Omega,\Omega)$ and $\gone(\Omega,\Omega^{\gpbl})$.
For $\scrB\sbst\Omega$,
the properties $\sfin(\Lambda,\scrB)$ are trivial
(see Appendix \ref{toostrong}).
Thus we cannot hope to have an equivalent item
``ONE does not have a winning strategy in $\gone(\Lambda,\scrB)$''
in the theorems dealing with these covers.
Fortunately, there exists an elegant technique to
deal with these cases without appealing to $\Lambda$.

\begin{lem}\label{GoneGameLemma}
Assume that $\scrB$ is countably thick.
For a space $X$, the following are equivalent.
\be
\i ONE does not have a winning strategy in $\gone(\Omega,\scrB)$,
\i ONE does not have a winning strategy in $\gone(\seq{\O_n},\scrB)$.
\ee
\end{lem}
\begin{proof}
We prove that $1\Impl 2$.
Assume $F$ is a strategy for ONE in $\gone(\seq{\O_n},\scrB)$
whose covers tree is $\{\cU_\sigma\}_{\sigma\in\N^*}$.
Define a strategy $\tilde F$ for ONE in $\gone(\Omega,\scrB)$ as follows:
The first move of ONE is $\tilde\cU_{\<\>}=\Union_{\sigma\in\N^*}\cU_\sigma$.
If TWO chooses $U_\sigma$, then ONE responds with
$$\tilde\cU_{\sigma}=\Union_{\eta\in\N^*}\cU_{\sigma\cat\eta},$$
etc. Now, a game lost by ONE according to the
strategy $\tilde F$ can be completed (by choosing the moves of TWO appropriately)
to a game lost
by ONE according to $F$ in $\gone(\seq{\O_n},\scrB)$.
As $\scrB$ is countably thick, this shows that $F$ is
not a winning strategy.
\end{proof}

\begin{thm}\label{GoneGameGeneral}
Fix $\scrB\in\{\Omega,\Omega^{\gpbl}\}$.
For each space $X$, the following are equivalent.
\be
\i $X$ satisfies $\sone(\Omega,\scrB)$,
\i ONE does not have a winning strategy in $\gone(\Omega,\scrB)$,
\i ONE does not have a winning strategy in $\gone(\seq{\O_n},\scrB)$.
\ee
\end{thm}
\begin{proof}
$2\Iff 3$ by Lemma \ref{GoneGameLemma}.

$1\Iff 2$:
For $\scrB=\Omega$ this is Theorem 2 of \cite{coc3}.
For $\scrB=\Omega^{\gpbl}$ this is Theorem 17 of \cite{coc7}.
\end{proof}

We do not know whether analogous game theoretic characterizations
can be given to the remaining few properties. The most interesting
problem seems to be the following.
\begin{prob}
Is it true that $X$ is a strong $\gamma$-set if, and only if,
ONE has no winning strategy in the game $\gone(\seq{\O_n},\Gamma)$?
\end{prob}

\section{The Borel case and the discrete case}\label{purecomb}
We need not stop in the case of open covers.
One important variant of open covers is that of
\emph{countable Borel covers}.
As in \cite{CBC}, one can translate all of the results
presented here to this case as well.
Another important variant is that of
\emph{arbitrary} countable covers of an uncountable cardinal $\kappa$.
Since these are exactly the
countable open covers of $\kappa$ with respect to the discrete topology on
$\kappa$, our results apply in this purely combinatorial case as well,
and we obtain new characterizations
of some well known combinatorial cardinal characteristics of
the continuum. For example:
\be
\i The \emph{unbounding number} $\b$ is equal to the minimal cardinal $\kappa$
such that ONE has a winning strategy in the game
$\gfin(\seq{\O_n},\O^{\gamma\gp})$ played on $\kappa$. 
\i The \emph{dominating number} $\d$ is equal to the minimal cardinal $\kappa$
such that ONE has a winning strategy in any (and both) of the games
$\gfin(\seq{\O_n},\Lambda)$ and 
$\gfin(\seq{\O_n},\O^{\omega\gp})$, 
played on $\kappa$.
\i The \emph{covering number for the meager ideal} $\cov(\M)$
is equal to the minimal cardinal $\kappa$
such that ONE has a winning strategy in any (and both) of the games:
$\gone(\seq{\O_n},\Lambda)$ and 
$\gone(\seq{\O_n},\Omega)$, 
played on $\kappa$.
\i The \emph{additivity number for the meager ideal} $\add(\M)$
is equal to the minimal cardinal $\kappa$
such that ONE has a winning strategy in the game
$\gone(\seq{\O_n},\O^{\gamma\gp})$,
played on $\kappa$.
\ee
All of these results follow easily from the equivalences
with the corresponding properties using the operators $\sone$ and $\sfin$,
together with the known critical cardinalities of these properties
-- see \cite{coc2}.

\appendix
\section{Too strong properties}\label{toostrong}

Assume that $\scrA$ and $\scrB$ are collections of covers of $X$.
We say that $X$ satisfies $\binom{\scrA}{\scrB}$ if
each element of $\scrA$ contains an element of $\scrB$ \cite{tautau}.
Clearly $\sfin(\scrA,\scrB)$ implies $\binom{\scrA}{\scrB}$.

\begin{prop}
Assume that $X$ is an infinite $T_1$ space, and fix $n\in\N$.
Then $X$ does not
satisfy any of the following properties:
\be
\i $\binom{\Lambda}{\O_2}$,
\i $\binom{\O_n}{\Lambda}$; and
\i $\binom{\O_n}{\O_{n+1}}$.
\ee
\end{prop}
\begin{proof}
(1) Fix a nonrepeating sequence $\seq{x_n}$ of elements of $X$,
and two distinct elements $a,b\in X$.
As $X$ is $T_1$, all singletons are closed subsets of $X$.
Then
$$\cU=\{X\sm\{x_{2n},a\} : n\in\N\}\cup\{X\sm\{x_{2n+1},b\} : n\in\N\}$$
is a large open cover of $X$, and for each $U\in\cU$, $a,b\not\sbst U$.

(2) and (3): Fix distinct elements $x_1,\dots,x_{n+1}\in X$.
Then
$$\cU=\{X\sm\{x_1\},\dots,X\sm\{x_{n+1}\}\}$$
is an open $n$-cover of $X$. But the $n+1$-element set
$\{x_1,\dots,x_{n+1}\}$ is not contained in any member of $\cU$.
As $\cU$ is a finite cover, it does not contain a large cover
either.
\end{proof}

Any nontrivial space has an open $\gamma$-cover:

\begin{lem}\label{alwaysgamma}
Assume that $X$ is an infinite $T_1$ space. Then $X$ has an open $\gamma$-cover.
\end{lem}
\begin{proof}
Fix a nonrepeating sequence $\seq{x_n}$ of element of $X$. Then
$$\cU=\{X\sm\{x_n\} : n\in\N\}$$
is an open $\gamma$-cover of $X$.
\end{proof}



\end{document}